\documentclass[a4paper,11pt]{article}
\usepackage[active]{srcltx}
\usepackage{amsmath}\usepackage{amssymb}\usepackage{amscd}

\usepackage[]{graphicx}

\newcommand{\Res}[1]{\lfloor #1\rfloor} 

\newcommand{\pr}[2]{\left\{#1,#2\right\}}

\newcounter{pic}

\newtheorem{defi}{Definition}

\newtheorem{prop}[defi]{Proposition}
\newtheorem{lemm}[defi]{Lemma}

\usepackage{empheq}

\begin{document}

\begin{center}

\bigskip
{\Large\bf Alternating subalgebras of Hecke algebras and alternating subgroups of braid groups}

\vspace{.6cm}

{\large {\bf O. V. Ogievetsky$^{\circ
}$\footnote{On leave of absence from P. N. Lebedev Physical Institute, Leninsky Pr. 53,
117924 Moscow, Russia} and L. Poulain d'Andecy$^{\circ}$\footnote{corresponding author (email adress:\ lpoulain@cpt.univ-mrs.fr)}}}

\vskip 0.6cm

$\circ\ $ Aix Marseille University, Center of Theoretical Physics, UMR 7332,\\
Luminy, 13288 Marseille, France

\end{center}

\vskip .6cm
\begin{abstract}\noindent
{}For a Coxeter system $(G,S)$ the multi-parametric alternating subalgebra $H^+(G)$ of the Hecke algebra and the alternating subgroup 
$\mathcal{B}^+(G)$ of the braid group are defined. 
Two presentations for $H^+(G)$ and $\mathcal{B}^+(G)$ are given; one generalizes the Bourbaki presentation for the alternating subgroups of Coxeter groups, another one uses generators related to edges of the Coxeter graph.
\end{abstract}

\section{\hspace{-0.55cm}.\hspace{0.55cm}Introduction}

Let $(G,S)$ be an arbitrary Coxeter system. The alternating subgroup $G^+$ is an index 2 subgroup of the Coxeter group $G$ (we will omit the reference to $S$, as in $G^+$ or, see below, in $H(G)$, etc). In \cite{Bou}, a presentation of the alternating group $G^+$ is suggested. 
In this ``Bourbaki" presentation,
one vertex of the Coxeter graph plays a particular role. In \cite{OPdA}, a different presentation has been given for alternating subgroups of Coxeter groups. In this presentation, the generators are related  
to the oriented edges of the Coxeter graph; for an irreducible Coxeter system $(G,S)$, no particular vertex is distinguished.

The Hecke algebra $H(G)$ associated to the Coxeter group $G$ is a flat deformation of the group ring of 
$G$. In \cite{Rat}, an analogue of the Hecke algebra is defined, in the one-parameter setting, for the alternating subgroups of the Coxeter groups. Here we extend this definition to the general multi-parameter situation and call the resulting algebra, denoted by $H^+(G)$, ``the alternating subalgebra of the Hecke algebra''.
The algebra $H^+(G)$ is an index 2 subalgebra (see Section \ref{sec-def} for precise definitions) of the Hecke algebra $H(G)$ and is a flat deformation of the group ring of $G^+$. The algebra $H^+(G)$ is among the deformations of the group ring of $G^+$ studied in \cite{ER2}.

In addition, associated to a Coxeter system $(G,S)$, there is a braid group $\mathcal{B}(G)$. Similarly to the alternating subgroup $G^+$ of $G$, we define an ``alternating subgroup" $\mathcal{B}^+(G)$ of the braid group $\mathcal{B}(G)$.
 
We give a presentation \`a la Bourbaki for the alternating subalgebra $H^+(G)$ of the Hecke algebra and
for the alternating subgroup $\mathcal{B}^+(G)$ of the braid group $\mathcal{B}(G)$. 
The Bourbaki presentation of $\mathcal{B}^+(G)$ (as well as the Bourbaki presentation of $G^+$) can be obtained by the Reidemeister--Schreier rewriting process \cite{Rei,Sch}; we present however a different proof.
Then we prove a presentation of $H^+(G)$ and  $\mathcal{B}^+(G)$ with generators related to the
oriented edges of the Coxeter graph.

One advantage of the presentation of $H^+(G)$
with generators related to edges of the Coxeter graph is that, passing from the defining relations for $G^+$ to the defining relations for $H^+(G)$, only the characteristic equations for the generators are 
deformed. This is similar to the situation of the Coxeter group $G$ and its Hecke algebra $H(G)$.

{}For type A, 
the presentations for the chain of algebras $H^+(A_n)$ and the chain of groups $\mathcal{B}^+(A_n)$
with generators related to edges of the Coxeter graph are local and stationary, in the sense of \cite{V}; for  
types B and D  
these presentations of the chains 
of algebras $H^+(B_n)$ and $H^+(D_n)$, and the chains of groups $\mathcal{B}^+(B_n)$ 
and $\mathcal{B}^+(D_n)$ are local and eventually stationary. This extends results obtained in \cite{OPdA}  for the alternating subgroups 
of Coxeter groups.

\vskip .1cm
The paper is organized as follows. In Section \ref{sec-def} we give the definition of the alternating subalgebra $H^+(G)$ of the Hecke algebra, as the even, for a certain grading, subalgebra
of the Hecke algebra $H(G)$. The Bourbaki presentation of $H^+(G)$ and the presentation of $H^+(G)$ with generators related to edges of the Coxeter graph are proved in Sections \ref{sec-bour} and \ref{sec-edge}. In Section \ref{sec-braid} we give analogues
of these two presentations for the alternating subgroups of the braid groups.
In Appendix, we obtain a recurrence relation and the generating function for the coefficients 
in the defining relations of the alternating subalgebras of the Hecke algebra.

\vskip .1cm
\noindent{\bf{Notation.}} 
Certain defining relations in this text involve a parameter $m\in \mathbb{Z}_{>0}\cup\infty$. It is understood that if $m=\infty$, the relation is absent.

\vskip .1cm\noindent
{}For any non-negative integer $k$, $\pr{a}{b}_k$ denotes the product $abab\dots$ with $k$ factors (by convention $\pr{a}{b}_0:=1$); for example $\pr{a}{b}_1:=a$, $\pr{a}{b}_2:=ab$ and $\pr{a}{b}_3:=aba$. We also set, for any non-negative integer $k$, $\pr{a}{b}_{-k}:=\pr{b}{a}_k$.

\section{\hspace{-0.55cm}.\hspace{0.55cm}Definition of the alternating subalgebra of the 
 Hecke algebra}\label{sec-def}

Let $(G,S)$ be a Coxeter system: $S$ is the set of generators, $S=\{ s_0,\dots,s_{n-1}\}$; 
the defining relations of the Coxeter group $G$ are encoded by a symmetric matrix $\mathfrak{m}$ with $m_{ii}=1$ and $2\leq m_{ij}\in
\mathbb{Z}_{>0}\cup\infty$ for $0\leq i<j\leq n-1$:
\begin{equation}\label{codecoxet} G=\langle\ s_0,\dots,s_{n-1}|\quad (s_is_j)^{m_{ij}}=1
\ \   \textrm{for $i,j=0,\dots,n-1$ such that $i\leq j$}\ \rangle.\end{equation} 
The sign character is the unique homomorphism $\epsilon : G\to \left\{-1,1\right\}$ such that $\epsilon(s_i)=-1$ for $i=0,\dots,n-1$. Its kernel $G^+:=\textrm{ker}(\epsilon)$ is called the alternating subgroup of $G$.

Recall that $s_i$ and $s_j$ are conjugate in the group $G$ iff there are some 
$i_1,\dots,i_r\in\{0,\dots,n-1\}$ such that $m_{ii_1},m_{i_1i_2},\dots,m_{i_rj}$ are odd, see, \emph{e.g.}, \cite{Hum}. 
 Let $(q_0,\dots,q_{n-1})$ be a set of indeterminates such that $q_i=q_j$ if $m_{ij}$ is odd. Let ${\cal{A}}$ be the ring of Laurent polynomials in 
$q_i$, $i=0,\dots ,n-1$, over $\mathbb{C}$.
The {\it Hecke algebra} $H(G)$ is the algebra over ${\cal{A}}$ generated by $g_0,\dots,g_{n-1}$ with the defining relations:
\begin{equation}\label{mult-Hecke1}
g_i^2=(q_i-q_i^{-1})g_i+1\ \   \textrm{for $i=0,\dots,n-1$,}\end{equation}
\begin{equation}\label{mult-Hecke2}
\pr{g_i}{g_j}_{m_{ij}}=\pr{g_j}{g_i}_{m_{ij}}\ \   \textrm{for $i,j=0,\dots,n-1$ such that $i<j$.}\end{equation}
The algebra $H(G)$ is a flat deformation of the group ring ${\cal{A}}G$: $H(G)$ has a basis whose elements are in one-to-one correspondence with the elements of $G$ (see chap. IV, sec. 1 exercise 23 in \cite{Bou}, \cite{C} and \cite{Eri} for 
different proofs).

The assignment $g_i\mapsto -g_i^{-1}$ extends to an involutive homomorphism  $\phi\colon H(G)\to H(G)$ of algebras, $\phi^2 =\text{id}$. Thus,
$H(G)=H^+(G)\oplus H^-(G)$, where $H^+(G)$ and $H^-(G)$ are eigenspaces of $\phi$ corresponding to eigenvalues $+1$ and
$-1$; the involution $\phi$ defines a $\mathbb{Z}_2$-grading on $H(G)$. The subalgebra $H^+(G)$ of even elements is called the 
{\it alternating subalgebra} of the Hecke algebra $H(G)$. 

Let ${\cal{B}}={\cal{B}}^+\oplus {\cal{B}}^-$ be a  $\mathbb{Z}_2$-graded associative algebra. Assume that ${\cal{B}}^-$ contains an invertible
element $f$. Then the left multiplication by $f$ gives an isomorphism of vector spaces ${\cal{B}}^+$ and ${\cal{B}}^-$. We then say that 
${\cal{B}}^+$ is a {\it subalgebra of index 2} of ${\cal{B}}$. 
Define the following elements of $H(G)$:
\[f_i:=\frac{1}{q_i+q_i^{-1}}(g_i+g_i^{-1})=\frac{2g_i-(q_i-q_i^{-1})}{q_i+q_i^{-1}},\quad\,i=0,\dots,n-1.\]
Since $f_i^2=1$ and $\phi (f_i)=-f_i$, $i=0,\dots,n-1$, $H^+(G)$ is a 
subalgebra of index 2 of $H(G)$; therefore, $H^+(G)$ is a flat 
deformation of the group ring ${\cal{A}}G^+$.

The elements $f_i$, $i=0,\dots,n-1$, form a generating set of $H(G)$. 
The algebra $H^+(G)$ is generated by the elements $f_if_j$, $i,j=0,\dots,n-1$ and $i\neq j$. 

Let $\beta_i:=
(q_i-q_i^{-1})/(q_i+q_i^{-1})$, $i=0,\dots,n-1$.  
Since the expression $\pr{g_i}{g_j}_{m_{ij}}-\pr{g_j}{g_i}_{m_{ij}}$  
is antisymmetric with respect to $i\leftrightarrow j$, 
the defining relations (\ref{mult-Hecke2}) of $H(G)$, in terms of the elements $f_i$, $i=0,\dots,n-1$, can be rewritten in the form
\begin{equation}\label{mult-Hecke-f'}
\sum\limits_{k=1}^{m_{ij}}a_k^{(m_{ij})}\bigl(\pr{f_i}{f_j}_k-\pr{f_j}{f_i}_k\bigr)+\sum\limits_{k=0}^{m_{ij}-1}b_k^{(m_{ij})}\bigl(\pr{f_i}{f_j}_k+\pr{f_j}{f_i}_k\bigr)=0\ .
\end{equation}
The leading coefficient $a_{m_{ij}}^{(m_{ij})}$ is non-zero and we normalize it to be $a_{m_{ij}}^{(m_{ij})}=1$. With this choice, 
$a_k^{(m_{ij})},b_k^{(m_{ij})} \in\mathbb{Z}[\beta_i,\beta_j]$ are polynomials in $\beta_i,\beta_j$ with integer coefficients;
$a_k^{(m_{ij})}$ is symmetric while $b_k^{(m_{ij})}$ is antisymmetric with respect to $\beta_i\leftrightarrow\beta_j$. 

\begin{lemm}\hspace{-.2cm}.\hspace{.2cm}
\label{lemm-akbk}
We have
\begin{equation}\label{lemm-a}
a_k^{(m_{ij})}=0\ \ \text{if\ \  $k\not\equiv m_{ij}\, ({\rm mod}\, 2)$},\end{equation}
\begin{equation}\label{lemm-a1}
b_k^{(m_{ij})}=0\ \ \text{\rm for any $k$}.\end{equation}
\end{lemm}
\emph{Proof.} The algebra 
${\cal{D}}_{ij}$ with the generators $g_i$ and $g_j$ and the defining relations (\ref{mult-Hecke1})--(\ref{mult-Hecke2}) is a 
flat deformation of the group ring of the dihedral group. The elements $1, \pr{g_i}{g_j}_k,\pr{g_j}{g_i}_k$, $k=1,\dots ,m_{ij}-1$,
and $\pr{g_i}{g_j}_{m_{ij}}$ form a basis $\mathfrak{B}$ of ${\cal{D}}_{ij}$. Denote by $\rho_{ij}$ the expression in the left hand side of (\ref{mult-Hecke-f'}). 
As $\phi$ is an automorphism of ${\cal{D}}_{ij}$, 
the relation $\phi (\rho_{ij})=0$ holds in the algebra ${\cal{D}}_{ij}$. Assume that there exists $k\not\equiv m_{ij} (\text{mod}\, 2)$ such that  
$a_k^{(m_{ij})}\neq 0$ or $b_k^{(m_{ij})}\neq 0$. Then $\rho_{ij}-(-1)^{m_{ij}}\phi (\rho_{ij})=0$ 
can be rewritten as a relation between the elements of $\mathfrak{B}$, a contradiction. Thus,
\begin{equation}\label{lemm-a'}
a_k^{(m_{ij})}=0\ \ 
\text{and}\ \ b_k^{(m_{ij})}=0\ \ 
\ \ \text{if\ \ $k\not\equiv m_{ij}\, ({\rm mod}\, 2)$}.\end{equation} 
If $m_{ij}$ is odd then $\beta_i=\beta_j$, so the expression $\pr{g_i}{g_j}_{m_{ij}}-\pr{g_j}{g_i}_{m_{ij}}$, rewritten in terms of the elements $f_i,f_j$, is antisymmetric with respect to $f_i\leftrightarrow f_j$;
thus $b_k^{(m_{ij})}$ vanish. Assume, for $m_{ij}$ even, that there exists an even $k$ 
such that $b_k^{(m_{ij})}\neq 0$. Then, taking into account (\ref{lemm-a'}), we can rewrite $\rho_{ij}+f_i\rho_{ij}f_i$ as 
a relation between the elements of $\mathfrak{B}$, a contradiction.\hfill$\square$ 

To conclude, we obtain the following set of defining relations of $H(G)$: 
\begin{equation}\label{mult-Hecke-f2}f _i^2=1\ \ \textrm{for $i=0,\dots,n-1$,}\end{equation}
and, multiplying the relation (\ref{mult-Hecke-f'}) by $(-1)^{m_{ij}}\pr{f_i}{f_j}_{m_{ij}}$, 
\begin{equation}\label{mult-Hecke-f22}
\sum\limits_{k=1}^{m_{ij}}a_k^{(m_{ij})}\Bigl((f_if_j)^{\frac{m_{ij}+k}{2}}-(f_if_j)^{\frac{m_{ij}-k}{2}}\Bigr)=0\ \ \textrm{for $i,j=0,\dots,n-1$ such that $i<j$}
\end{equation}
with the restriction (\ref{lemm-a}). 

{}For the types A and B, the algebra $H^+(G)$ was introduced in \cite{Mit,Mit2}. 
{}For equal parameters, $q_i=q$, the coefficients $a_k^{(m_{ij})}$ have been calculated in \cite{Rat}. 
We study the coefficients $a_k^{(m_{ij})}$ with arbitrary $q_i$ and $q_j$ in the Appendix. 

\section{\hspace{-0.55cm}.\hspace{0.55cm}Bourbaki presentation}\label{sec-bour}

Let $(G,S)$ be a Coxeter system with the Coxeter matrix $\mathfrak{m}$. We first recall the Bourbaki presentation of $G^+$ suggested in \cite{Bou}, chap. IV, sec. 1, exercise 9 (see \cite{BRR} for a proof). The alternating group $G^+$ is isomorphic to the group generated by $R_1,\dots,R_{n-1}$ with the defining relations:
\begin{equation}\label{codebourb}\left\{\begin{array}{ll}
R_i^{m_{0i}}=1 & \textrm{for $i=1,\dots,n-1$,}\\[.2em]
(R_i^{-1}R_j)^{m_{ij}}=1 & \textrm{for $i,j=1,\dots,n-1$ such that $i<j$.}\end{array}\right.\end{equation}
The isomorphism with $G^+$ is given by $R_i\mapsto s_0s_i$ for $i=1,\dots,n-1$. 
The Bourbaki presentation of the alternating group $G^+$ depends on the choice of a generator carrying the subscript 0. 

We prove in this Section a presentation of $H^+(G)$ similar to the Bourbaki presentation (\ref{codebourb}) of $G^+$.

\begin{prop}\hspace{-.2cm}.\hspace{.2cm}
\label{pres-bour}
For a Coxeter system $(G,S)$ with the Coxeter matrix $\mathfrak{m}$, the alternating subalgebra $H^+(G)$ of the Hecke algebra is isomorphic to the algebra $\mathfrak{A}$
with the generators $Y_1^{\pm1},\dots,Y_{n-1}^{\pm1}$ and the defining relations
\begin{equation}\label{rel-bour}\left\{\begin{array}{ll}
\sum\limits_{k=1}^{m_{0i}}a_k^{(m_{0i})}\Bigl(Y_i^{\frac{m_{0i}+k}{2}}-Y_i^{\frac{m_{0i}-k}{2}}\Bigr)=0 & \textrm{for $i=1,\dots,n-1$,}\\[.2em]
\sum\limits_{k=1}^{m_{ij}}a_k^{(m_{ij})}\Bigl((Y_i^{-1}Y_j)^{\frac{m_{ij}+k}{2}}-(Y_i^{-1}Y_j)^{\frac{m_{ij}-k}{2}}\Bigr)=0 & \textrm{for $i,j=1,\dots,n-1$ such that $i<j$.}\end{array}\right.\end{equation}
\end{prop}
\emph{Proof.}
Define a map $\psi$ from the set of generators $\{Y_1,\dots,Y_{n-1}\}$ to the algebra $H^+(G)$ by
\[Y_i\mapsto f_0f_i\ \text{for $i=1,\dots,n-1$.}\]
Due to the relations (\ref{mult-Hecke-f2})--(\ref{mult-Hecke-f22}), this map extends to a (surjective) homomorphism, which we denote again by 
$\psi$, from the algebra $\mathfrak{A}$ to $H^+(G)$. We shall prove that $\psi$ is an isomorphism.

The left hand side of the first relation in (\ref{rel-bour}) is invariant under the following sequence of operations: replace $Y_i$ by $Y_i^{-1}$ and then multiply by $-Y_i^{m_{0i}}$. One can verify directly that the left hand side of the second relation in (\ref{rel-bour}) is invariant under the following sequence of operations: replace $Y_i$ by $Y_i^{-1}$, $Y_j$ by $Y_j^{-1}$, then multiply from the left by $-Y_j^{-1}(Y_jY_i^{-1})^{m_{ij}}$ and from the right by $Y_j$. Therefore the map defined by $\omega\,:\,Y_i\mapsto Y_i^{-1}$ extends to an involutive automorphism of the algebra $\mathfrak{A}$.

With the help of $\omega$, we define the cross-product $\tilde{\mathfrak{A}}$ of the algebra $\mathfrak{A}$ with the cyclic group $C_2$ with two elements.  
Let $f$ be the generator of the group $C_2$. As a vector space, $\tilde{\mathfrak{A}}$ is isomorphic to $\mathfrak{A}\otimes\mathcal{A}C_2$. The generators of $\tilde{\mathfrak{A}}$ are the elements $Y_1^{\pm1},\dots,Y_{n-1}^{\pm1}$ with the defining relations (\ref{rel-bour}), and in addition the generator $f$ with the defining relations $f^2=1$ and $fY_i=Y_i^{-1}f$, $i=1,\dots,n-1$. The map
\[f\mapsto f_0\ \ \ \text{and}\ \ \ Y_i\mapsto f_0f_i\ \ \text{for $i=1,\dots,n-1$},\]
extends to a morphism of algebras $\psi_1\,:\,\tilde{\mathfrak{A}}\to H(G)$. The verification is straightforward (use 
(\ref{mult-Hecke-f2})-(\ref{mult-Hecke-f22})). On the other hand, one directly verifies that the map
\[f_0\mapsto f\ \ \ \text{and}\ \ \ f_i\mapsto f Y_i\ \ \text{for $i=1,\dots,n-1$},\]
extends to a morphism of algebras $\psi_2\,:\,H(G)\to\tilde{\mathfrak{A}}$. Moreover, the morphisms $\psi_1$ and $\psi_2$ are mutually inverse. The restriction of $\psi_2$ to $H^+(G)$ is the morphism inverse to $\psi$.\hfill$\square$

\section{\hspace{-0.55cm}.\hspace{0.55cm}Presentation using edges of the Coxeter graph}\label{sec-edge}

Let $(G,S)$ be a Coxeter system with the Coxeter matrix $\mathfrak{m}$. We first recall the presentation given in \cite{OPdA} of $G^+$; it uses edges of the Coxeter graph ${\cal{G}}$ of $(G,S)$.

Recall that vertices,  indexed by the subscripts $0,1,\dots ,n-1$, of the Coxeter graph ${\cal{G}}$  
are in one-to-one correspondence with the generators $s_0,\dots,s_{n-1}$ of $G$; 
vertices $i$ and $j$ are connected if and only if $m_{ij}\geq3$ and then the edge between $i$ and $j$ is labeled by the number $m_{ij}$.
In the sequel the edge between vertices $i$ and $j$ is denoted by $(ij)$. 

If ${\cal{G}}$ is not connected, let ${\cal{G}}_1,{\cal{G}}_2,\dots ,{\cal{G}}_m$ be its connected components.
We choose an arbitrary vertex $i_a$ of ${\cal{G}}_a$ for $a=1,\dots,m$;  we add an edge between $i_l$ and 
$i_{l+1}$ for $l=1,\dots,m-1$ and label it by the number $2$. 
The obtained connected graph  ${\cal{G}}^c$ we call a \emph{connected extension} of the Coxeter graph $\mathcal{G}$.

The presentation in \cite{OPdA} uses an orientation - chosen arbitrarily - of edges 
of the connected extension ${\cal{G}}^c$ of the Coxeter graph. 
For concreteness, if there is an edge between $i$ and $j$ with $i<j$, we orient it from $i$ to $j$. 
We associate a generator $r_{ij}$ to each oriented edge of 
${\cal{G}}^c$. For a generator $r_{ij}$ we denote by $r_{ji}$ the inverse, $r_{ji}:=r_{ij}^{-1}$.
\begin{defi}\hspace{-.2cm}.\hspace{.2cm}
\label{def-dis-edge2}
Two edges $(ij)$ and $(kl)$ of ${\cal{G}}^c$ are said to be not connected if $\left\{i,j\right\}\cap\left\{l,m\right\}=\varnothing$ and
there is no edge connecting any of the vertices $\left\{i,j\right\}$ with any of the vertices $\left\{l,m\right\}$.
\end{defi}
The alternating group $G^+$ is isomorphic \cite{OPdA} to the group with the generators $r_{ij}$ and the defining relations
\begin{equation}\label{code-edge}\left\{\begin{array}{lll}
(r_{ij})^{m_{ij}}=1 && \textrm{for all generators $r_{ij}$,}\\[.2em]
r_{ii_1}r_{i_1i_2}\dots r_{i_ai}=1 && \textrm{for cycles with edges  $(ii_1),(i_1i_2)\dots,(i_ai)$}, 
\\[.2em]
(r_{ij}r_{jk})^2=1 && \textrm{for $r_{ij},r_{jk}$ such that $i<k$ and $m_{ik}=2$,}\\[.2em]
(r_{ij}r_{jk}r_{kl})^2=1 && \textrm{for $r_{ij},r_{jk},r_{kl}$ such that $i<l$ and $m_{il}=2$,}\\[.2em]
r_{ij}r_{lm}=r_{lm}r_{ij} && \textrm{if $(ij)$ and $(lm)$ are not connected.}
\end{array}\right.\end{equation}
The isomorphism with $G^+$ is given by $r_{ij}\mapsto s_is_j$ for all generators $r_{ij}$.

We generalize this presentation to a presentation of $H^+(G)$. Associate an element $y_{ij}$ to each generator $r_{ij}$ of $G^+$, and set, for all $y_{ij}$, $y_{ji}:=y_{ij}^{-1}$. 

\begin{prop}\hspace{-.2cm}.\hspace{.2cm}
\label{pres-edge}
The alternating subalgebra $H^+(G)$ of the Hecke algebra is isomorphic to the algebra $\mathfrak{Y}$
with the generators $y_{ij}$ and the defining relations
\begin{empheq}[left=\empheqlbrace]{alignat=1}\label{rel-edge-a}&\sum\limits_{k=1}^{m_{ij}}a_k^{(m_{ij})}(y_{ij}^{\frac{m_{ij}+k}{2}}-y_{ij}^{\frac{m_{ij}-k}{2}})=0 \hspace{1.0cm} \textrm{for all generators $y_{ij}$,}\\[.1em]
\label{rel-edge-b}&y_{ii_1}y_{i_1i_2}\dots y_{i_ai}=1 \hspace{2.9cm} \textrm{for cycles with edges  $(ii_1),(i_1i_2)\dots,(i_ai)$}, 
\\[.1em]
\label{rel-edge-c}&(y_{ij}y_{jk})^2=1 \hspace{3.9cm} \textrm{for $y_{ij},y_{jk}$ such that $i<k$ and $m_{ik}=2$,}\\[.1em]
\label{rel-edge-d}&(y_{ij}y_{jk}y_{kl})^2=1 \hspace{3.4cm} \textrm{for $y_{ij},y_{jk},y_{kl}$ such that $i<l$ and $m_{il}=2$,}\\[.1em]
\label{rel-edge-e}&y_{ij}y_{lm}=y_{lm}y_{ij} \hspace{3.45cm} \textrm{if $(ij)$ and $(lm)$ are not connected.}\end{empheq}
\end{prop}
Let $\mathfrak{c}_1,\dots,\mathfrak{c}_{\mathfrak{l}}$ be a set of generators of the fundamental group of ${\cal{G}}^c$. 
In the set of the defining relations it is sufficient to impose the relation (\ref{rel-edge-b}) for the cycles $\mathfrak{c}_{\mathfrak{a}}$, $\mathfrak{a}=1,\dots,\mathfrak{l}$.

\noindent\emph{Proof of the Proposition.} Notice that if $m_{0i}=2$ then the first relation in (\ref{rel-bour}) reduces to $Y_i^2=1$ and also that, if $m_{ij}=2$ then the second relation in (\ref{rel-bour}) reduces to $(Y_i^{-1}Y_j)^2=1$. Due to this fact, the proof is very similar to the proof in \cite{OPdA} in the classical situation (that is, for the presentation (\ref{code-edge}) of $G^+$). So we only sketch it.  
The following map
\begin{equation}\label{real}y_{ij}\mapsto \left\{\begin{array}{ll}Y_i^{-1}Y_j & \textrm{if $i\neq0$,}\\[.2em] 
Y_j & \textrm{if $i=0$,}\end{array}\right.\end{equation}
extends to an algebra homomorphism $\Phi\,:\,\mathfrak{Y}\to H^+(G)$.

Define now the map $\Psi$ from the set of generators $\{Y_1,\dots,Y_{n-1}\}$ of $H^+(G)$ to $\mathfrak{Y}$ by
\begin{equation}\label{inv-real}\Psi\ :\ Y_i\mapsto{\grave{Y}}_i:=y_{0i_1}y_{i_1i_2}\dots y_{i_ki}\ \ \textrm{for all $i=1,\dots,n-1$,}\end{equation}
where $(0,i_1,i_2,\dots,i_k,i)$ is an arbitrary path from the vertex $0$ to the vertex $i$ in ${\cal{G}}^c$. 
The map $\Psi$ is well-defined since the element $\grave{Y}_i$ does not depend on the chosen path, due to the relation (\ref{rel-edge-b}). The map $\Psi$ extends to an algebra homomorphism from $H^+(G)$ to $\mathfrak{Y}$, which we still denote by $\Psi$. Moreover $\Psi$ and $\Phi$ are mutually inverse.\hfill$\square$

\vskip .1cm\noindent
\textbf{Remark.} The defining relations (\ref{rel-edge-a})--(\ref{rel-edge-e}) of the 
algebra $H^+(G)$ are deformations of the defining relations (\ref{code-edge}) of the 
group $G^+$. Only the characteristic equation for the generators is deformed. This is similar to the Hecke algebra situation (passing from the relations (\ref{codecoxet}) to the relations (\ref{mult-Hecke1})--(\ref{mult-Hecke2}), only the characteristic equation for the generators is deformed). This phenomenon does not appear in the deformation of
the Bourbaki presentation (\ref{codebourb}) of $G^+$ to the Bourbaki presentation (\ref{rel-bour}) of $H^+(G)$.

\section{\hspace{-0.55cm}.\hspace{0.55cm}Alternating subgroups of braid groups}\label{sec-braid}

\subsection{Definition}

Let $(G,S)$ be a Coxeter system with the Coxeter matrix $\mathfrak{m}$.
The braid group $\mathcal{B}(G)$ is the group generated by $g_0,\dots,g_{n-1}$ with the defining relations:
\begin{equation}\label{braid}
\pr{g_i}{g_j}_{m_{ij}}=\pr{g_j}{g_i}_{m_{ij}}\ \ \ \textrm{for $i,j=0,\dots,n-1$ such that $i<j$.}\end{equation}
Extend the sign character to the group $\mathcal{B}(G)$, that is, define the homomorphism $\epsilon : \mathcal{B}(G)\to \left\{-1,1\right\}$ by $\epsilon(g_i)=-1$ for $i=0,\dots,n-1$. The kernel, $\mathcal{B}^+(G):=\textrm{ker}(\epsilon)$, we call the alternating subgroup of the braid group $\mathcal{B}(G)$.
The group $\mathcal{B}^+(G)$ is generated by the elements $g_ig_j$ (and their inverses), 
$i,j=0,\dots,n-1$.

\vskip .1cm\noindent
\textbf{Remark.} 
There is a natural $\mathbb{Z}_2$-grading of the group ring of $\mathcal{B}(G)$ defined by $\epsilon$.  
Let $\pi$ be the natural surjection of the group ring of $\mathcal{B}(G)$ to the Hecke algebra $H(G)$ (the quotient  
by the relation (\ref{mult-Hecke1})). Recall the grading on $H(G)$ defined by the involution $\phi$. It should be noted that 
$\phi\pi\neq \pi\epsilon$, the image of
$\mathcal{B}^+(G)$ under $\pi$ does not belong to $H^+(G)$; in other words, the grading on $H(G)$ is not induced by $\pi$ from the grading on the group ring of $\mathcal{B}(G)$. 

\subsection{Bourbaki presentation for alternating subgroups of braid groups}

We extend the Bourbaki presentations (\ref{codebourb}) and (\ref{rel-bour}) of the group $G^+$ and the algebra $H^+(G)$ to the group $\mathcal{B}^+(G)$. The presentation depends on a choice of a generator $g_0$, carrying the subscript $0$, among the generators of $\mathcal{B}(G)$.
\begin{prop}\hspace{-.2cm}.\hspace{.2cm}
\label{pres-braid}
For a Coxeter system $(G,S)$ with the Coxeter matrix $\mathfrak{m}$, the alternating subgroup $\mathcal{B}^+(G)$ of the braid group is isomorphic to the group $B$
with the generators $\mathsf{R}_0,\dots,\mathsf{R}_{n-1}$ and $\mathsf{R}'_0,\dots,\mathsf{R}'_{n-1}$ and the defining relations
\begin{equation}\label{rel-braid}\left\{\begin{array}{ll}
\mathsf{R}'_0=1\,,\\[.3em]
{\{\mathsf{R}'_i,\mathsf{R}_j\}}_{m_{ij}}={\{\mathsf{R}'_j,\mathsf{R}_i\}}_{m_{ij}} & \textrm{for $i,j=0,\dots,n-1$ such that $i<j$,}\\[.3em]
{\{\mathsf{R}_i,\mathsf{R}'_j\}}_{m_{ij}}={\{\mathsf{R}_j,\mathsf{R}'_i\}}_{m_{ij}} & \textrm{for $i,j=0,\dots,n-1$ such that $i<j$.}\end{array}\right.\end{equation}
\end{prop}
\emph{Proof.} 
Define a map $\psi$ from the set of generators of $B$ to $\mathcal{B}^+(G)$ by
\begin{equation}\label{iso-braid}\mathsf{R}_i\mapsto g_0g_i,\ i=0,\dots n-1,\ \ \text{and}\ \ \mathsf{R}'_i\mapsto g_ig_0^{-1},\ i=0,\dots n-1.\end{equation}
One directly verifies that $\psi$ extends to a homomorphism, which we still denote by $\psi$, from $B$ to $\mathcal{B}^+(G)$. Moreover, the homomorphism $\psi$ is surjective. Indeed, for $i,j=0,\dots,n-1$, we have $g_ig_j=\psi(\mathsf{R}'_i\mathsf{R}_j)$. 
We shall prove that $\psi$ is actually an isomorphism.

Define a map $\omega$ from the the set of generators of $B$ to $B$ by
\[\mathsf{R}_i\mapsto \mathsf{R}'_i\mathsf{R}_0,\ i=0,\dots n-1,\ \ \text{and}\ \ \mathsf{R}'_i\mapsto \mathsf{R}_0^{-1}\mathsf{R}_i,\ i=0,\dots n-1.\]
It is straightforward to verify that $\omega$ defines an automorphism of $B$. Moreover, the automorphism $\omega^2$ is inner: 
for any $x\in B$, we have $\omega^2(x)=\mathsf{R}_0^{-1}x\mathsf{R}_0$. 

The automorphism $\omega$ generates the action of the infinite cyclic group $C$ on $B$.
Let $B\rtimes C$ be the corresponding semidirect product;
the group $B\rtimes C$ is generated by the generators of $B$ and an element $g$, and we add to the defining relations of $B$ the relation 
$gxg^{-1}=\omega(x)$ for each generator $x$ of $B$. 

Now let $Q$ be the quotient of the group $B\rtimes C$ by the relation $g^2=\mathsf{R}_0^{-1}$. The following map
\[g\mapsto g_0^{-1},\ \ \mathsf{R}_i\mapsto g_0g_i,\ i=0,\dots,n-1,\ \ \text{and}\ \ \mathsf{R}'_i\mapsto g_ig_0^{-1},\ i=0,\dots,n-1,\]
extends to a homomorphism $\psi_1$ from $Q$ to $\mathcal{B}(G)$. The verification is the same as for the map $\psi$, given by (\ref{iso-braid}), with, in addition, the verification of the relations in $Q$ concerning the generator $g$; these are satisfied by construction.

The following map
\[g_i\mapsto g^{-1} \mathsf{R}_i,\ i=0,\dots,n-1,\]
extends to a homomorphism $\psi_2$ from $\mathcal{B}(G)$ to $Q$. We omit the straightforward calculations here.

Moreover, the morphisms $\psi_1$ and $\psi_2$ are mutually inverse. The restriction of $\psi_2$ to $\mathcal{B}^+(G)$ is the inverse of 
$\psi$; thus, the homomorphism $\psi$, given by (\ref{iso-braid}), is the required isomorphism between $B$ and $\mathcal{B}^+(G)$.\hfill$\square$

\vskip .1cm \noindent
\textbf{Remarks.} \textbf{(i)} Denote by $\tau$ the standard anti-automorphism of the group $\mathcal{B}(G)$, sending $g_i$ to $g_i$, $i=0,\dots,n-1$.  The action of $\tau$ on the generators of the Bourbaki presentation is given by
\[\mathsf{R}_i\to\mathsf{R}'_i\mathsf{R}_0\ \ \text{and}\ \ \mathsf{R}'_i\to\mathsf{R}_0^{-1}\mathsf{R}_i,\ \ \ i=0,\dots,n-1.\]

\textbf{(ii)} The Coxeter group $G$ is the quotient of the braid group $\mathcal{B}(G)$ by the relations $g_i^2=1$, $i=0,\dots,n-1$. In the alternating setting, we have a similar result; the alternating subgroup $G^+$ of the Coxeter group $G$ is the quotient of the group $\mathcal{B}^+(G)$ by the relations $\mathsf{R}_0=1$ and $\mathsf{R}'_i=\mathsf{R}_i^{-1}$, $i=1,\dots,n-1$. Indeed, in this quotient, the relations (\ref{rel-braid}) reduce to
\[\left\{\begin{array}{ll}
\mathsf{R}_j^{m_{0j}}=1 & \text{for $j=1,\dots,n-1$,}\\[0.3em]
(\mathsf{R}_i^{-1}\mathsf{R}_j)^{m_{ij}}=1 & \text{for $i,j=1,\dots,n-1$ such that $i<j$.}
\end{array}\right.\]
These are the defining relations of the Bourbaki presentation of $G^+$, see (\ref{codebourb}).

\vskip .1cm
\textbf{(iii)}
The Reidemeister--Schreier rewriting process \cite{Rei,Sch} (see \emph{e.g.} \cite{LS} for a more recent exposition),
allows to find a presentation of a subgroup $H$ of a group $G$, given a presentation of $G$ and a suitable information about $H$. We apply this process to the subgroup $\mathcal{B}^+(G)$ of $\mathcal{B}(G)$.
Decompose $\mathcal{B}(G)$ into the disjoint union of its right cosets with respect to $\mathcal{B}^+(G)$,  $\mathcal{B}(G)=\mathcal{B}^+(G)\cup \mathcal{B}^+(G)g_0$. For any $a\in\mathcal{B}(G)$, define $\overline{a}\in\{1,g_0\}$ by $\mathcal{B}^+(G)a=\mathcal{B}^+(G)\overline{a}$
and let ${\sf S}:=\{ g_0,g_1,\dots ,g_{n-1}\}$.
The Reidemeister--Schreier rewriting process asserts that $\mathcal{B}^+(G)$ is isomorphic to the group with a set of generators $\mathfrak{S}$ and a set of defining relations $\mathfrak{D}$ defined
as follows.
 
 \vskip .1cm\noindent
 $\bullet$ Elements of $\mathfrak{S}$ are in one-to-one correspondence with elements $\gamma(a,g):=ag(\overline{ag})^{-1}$, $g\in{\sf S}$
 and $a\in\{1,g_0\}$, such that $ag(\overline{ag})^{-1}\neq 1$; we obtain generators $\mathfrak{R}_i\in\mathfrak{S}$, $i=0,\dots,n-1$, corresponding to $g_0g_i(\overline{g_0g_i})^{-1}=g_0g_i$,
 and generators $\mathfrak{R}'_i\in\mathfrak{S}$, $i=1,\dots,n-1$, corresponding to $g_i(\overline{g_i})^{-1}=g_ig_0^{-1}$.
 Define for convenience $\mathfrak{R}'_0:=1$.

 \vskip .1cm\noindent
 $\bullet$ 
 For a word $w\! =\! g_{i_1}\dots g_{i_k}$ in the the alphabet $\sf{S}$  
let $\pi(w)\!  :=\! \hat{\gamma}(1,g_{i_1})\hat{\gamma}(\overline{g_{i_1}},g_{i_2})\dots\hat{\gamma}(\overline{g_{i_1}\dots g_{i_{k-1}}},g_{i_k})$ where $\hat{\gamma}(1,g_0):=\mathfrak{R}'_0=1$ and, otherwise, $\hat{\gamma}(a,g)$ is the generator corresponding to  $\gamma(a,g)$. 
The relations $\pi(a\pr{g_i}{g_j}_{m_{ij}})\!  =\! \pi(a\pr{g_j}{g_i}_{m_{ij}})$, $a\in\{1,g_0\}$ and $i,j=0,\dots,n-1$, $i<j$, form
the set $\mathfrak{D}$. It is straightforward to see that the relations in $\mathfrak{D}$  
are ${\{\mathfrak{R}'_i,\mathfrak{R}_j\}}_{m_{ij}}={\{\mathfrak{R}'_j,\mathfrak{R}_i\}}_{m_{ij}}$ and
${\{\mathfrak{R}_i,\mathfrak{R}'_j\}}_{m_{ij}}={\{\mathfrak{R}_j,\mathfrak{R}'_i\}}_{m_{ij}}$, $i,j=0,\dots,n-1$, $i<j$.

\vskip .1cm\noindent
Thus, the presentation of 
the Proposition \ref{pres-braid} coincides with the one 
obtained by the Reidemeister--Schreier rewriting process (with $\{1,g_0\}$ as the ``Schreier transversal")
for the subgroup $\mathcal{B}^+(G)$ of $\mathcal{B}(G)$.
{}For the alternating subgroup $G^+$ of the Coxeter group $G$ the Reidemeister--Schreier rewriting process leads to the presentation (\ref{codebourb}).

\vskip .1cm
\textbf{(iv)} The set of generators in the presentation of $\mathcal{B}^+(G)$ given by the Proposition \ref{pres-braid} is, in general, not minimal. For example, if there is some $j\in\{1,\dots,n-1\}$ such that $m_{0j}=2$, then the second and third relations in (\ref{rel-braid}) (for $i=0$) imply that $\mathsf{R}'_j=\mathsf{R}_j\mathsf{R}_0^{-1}=\mathsf{R}_0^{-1}\mathsf{R}_j$. Furthermore, if there is some $j\in\{1,\dots,n-1\}$ such that $m_{0j}$ is odd, then the second and third relations in (\ref{rel-braid}) (for $i=0$) imply that $\mathsf{R}_j^{\frac{m_{0j}-1}{2}}=(\mathsf{R}'_j\mathsf{R}_0)^{\frac{m_{0j}-1}{2}}\mathsf{R}'_j$ and $\mathsf{R}_0(\mathsf{R}'_j\mathsf{R}_0)^{\frac{m_{0j}-1}{2}}=\mathsf{R}_j^{\frac{m_{0j}+1}{2}}$. Thus, we have that
\[\mathsf{R}'_j=\mathsf{R}_j^{-\frac{m_{0j}+1}{2}}\mathsf{R}_0\mathsf{R}_j^{\frac{m_{0j}-1}{2}}\ \ \text{for $j=1,\dots,n-1$ such that $m_{0j}$ is odd.}\]
Nevertheless, in general, both sets, $\mathsf{R}_i$  and  $\mathsf{R}'_i$, of generators are needed. Consider, for example, the braid group 
$U$ generated by $g_0$ and $g_1$ and the defining relations $g_0g_1g_0g_1=g_1g_0g_1g_0$ (that is, $\pr{g_0}{g_1}_4=\pr{g_1}{g_0}_4$). In this case, the generators for the alternating subgroup suggested by the Proposition \ref{pres-braid} are $g_0^2$, $g_0g_1$ and $g_1g_0^{-1}$. The element $g_1g_0^{-1}$ does not belong to the subgroup generated by $g_0^2$ and $g_0g_1$. Indeed, let $\bar{U}$ be the quotient of $U$ by the normal subgroup generated by $g_0^2$ and let $\bar{g}_i$ be the images of $g_i$ in $\bar{U}$.
It is known that $\bar{U}$ is isomorphic to $S_2\ltimes C^2$, where $S_2$ is the symmetric group on 2 elements and $C$ is the infinite cyclic group; $S_2$ acts on $C^2$ by permuting the two copies of $C$. 
Suppose that there exists an integer $x$ such that $\bar{g}_1\bar{g}_0=(\bar{g}_0\bar{g}_1)^x$. If $x=1$, that is $\bar{g}_1\bar{g}_0=
\bar{g}_0\bar{g}_1$, then the group $\bar{U}$ would be isomorphic to $S_2\times C$, which is impossible. Assume that $x\neq 1$. We have 
$(\bar{g}_1\bar{g}_0)^2=(\bar{g}_0\bar{g}_1)^{2x}$ which, together with the defining relation $(\bar{g}_1\bar{g}_0)^2=
(\bar{g}_0\bar{g}_1)^2$ leads to $(\bar{g}_0\bar{g}_1)^{2(x-1)}=1$, contradicting to the fact that $\bar{U}$ is infinite.
A similar calculation shows that the element $g_0g_1\in U$ does not belong to the subgroup generated by $g_0^2$ and $g_1g_0^{-1}$. 

\vskip .1cm
\textbf{(v)} The elements $\mathsf{R}_i$, $i=0,\dots,n-1$, generate the alternating subgroup of the braid group of any simply-laced type, see
the remark \textbf{(v)}. We give the presentation for  the alternating subgroup of the braid group of type A.
Label the generators of the braid group $\mathcal{B}(A_n)$ in the standard way; that is, $\mathcal{B}(A_n)$ is generated by $g_0,g_1,\dots,g_{n-1}$ with the defining relations:
\begin{equation}\label{braid-A}
\left\{\begin{array}{ll}
g_ig_{i+1}g_i=g_{i+1}g_ig_{i+1} & \text{for $i=0,\dots,n-2$,}\\[0.3em]
g_ig_j=g_jg_i & \text{for $i,j=0,\dots,n-1$ such that $|i-j|>1$.}
\end{array}\right.
\end{equation}
The group $\mathcal{B}^+(A_n)$ is isomorphic to the group generated by $\mathsf{R}_0,\mathsf{R}_1,\dots,\mathsf{R}_{n-1}$ with the defining relations:
\begin{equation}\label{braid-A+}
\left\{\begin{array}{ll}
\mathsf{R}_0\mathsf{R}_1\mathsf{R}_0=\mathsf{R}_1^2\mathsf{R}_0^{-1}\mathsf{R}_1^2,\\[0.3em]
\mathsf{R}_0\mathsf{R}_j=\mathsf{R}_j\mathsf{R}_0 & \text{for $j=2,\dots,n-1$,}\\[0.3em]
\mathsf{R}_2\mathsf{R}_1\mathsf{R}_2=\mathsf{R}_1^2\mathsf{R}_0^{-1}\mathsf{R}_2
\mathsf{R}_0^{-1}\mathsf{R}_1^2,\\[0.3em]
\mathsf{R}_2\mathsf{R}_1^2\mathsf{R}_2=\mathsf{R}_0\mathsf{R}_1\mathsf{R}_2\mathsf{R}_0^{-1}
\mathsf{R}_1\mathsf{R}_0,\\[0.3em]
\mathsf{R}_1^2\mathsf{R}_j=\mathsf{R}_j\mathsf{R}_1\mathsf{R}_0,\ \ \mathsf{R}_j\mathsf{R}_1^2=\mathsf{R}_0\mathsf{R}_1\mathsf{R}_j & \text{for $j=3,\dots,n-1$,}\\[0.3em]
\mathsf{R}_i\mathsf{R}_{i+1}\mathsf{R}_i=\mathsf{R}_{i+1}\mathsf{R}_i\mathsf{R}_{i+1} & \text{for $i=2,\dots,n-2$,}\\[0.3em]
\mathsf{R}_i\mathsf{R}_j=\mathsf{R}_j\mathsf{R}_i & \text{for $i,j=2,\dots,n-1$ such that $|i-j|>1$.}\\[0.3em]
\end{array}\right.
\end{equation}
The verification that this presentation is equivalent to (\ref{rel-braid}) for the type A is straightforward once one notices that here we have $\mathsf{R}'_1=\mathsf{R}_0^{-1}\mathsf{R}_1^2\mathsf{R}_0^{-1}$ and $\mathsf{R}'_j=\mathsf{R}_j\mathsf{R}_0^{-1}=\mathsf{R}_0^{-1}\mathsf{R}_j$ for $j=2,\dots,n-1$.

It is interesting to note that in terms of generators $\bar{\mathsf{R}}_0:=\mathsf{R}_0^{-1}$ and $\mathsf{R}_i$, $i=1,\dots , n-1$, one can rewrite all relations (\ref{braid-A+}) without inverses of generators and define thus a monoid of positive elements.

\subsection{Presentation using edges of the Coxeter graph of alternating subgroups of braid groups}

The  group $\mathcal{B}^+(G)$ admits a presentation similar to the presentations of the group $G^+$ and the algebra $H^+(G)$, see  (\ref{code-edge}) and the Proposition \ref{pres-edge}. Associate, as in Section \ref{sec-edge}, a generator $\mathsf{r}_{ij}$ to any oriented edge (the edges are oriented from $i$ to $j$ if $i<j$) of the graph $\mathcal{G}^c$. Set $\mathsf{r}_{ji}:=\mathsf{r}_{ij}^{-1}$ for all generators $\mathsf{r}_{ij}$.

\begin{prop}\hspace{-.2cm}.\hspace{.2cm}
\label{prop-edge-br}
For a Coxeter system $(G,S)$ with the Coxeter matrix $\mathfrak{m}$, the alternating subgroup $\mathcal{B}^+(G)$ of the braid group is isomorphic to the group $\beth$ generated by the elements $\mathsf{r}_{ij}$ and $\mathsf{t}_0,\dots,\mathsf{t}_{n-1}$ with the defining relations
\begin{equation}\label{rel-braid2}\left\{\begin{array}{lll}
\mathsf{r}_{ii_1}\mathsf{r}_{i_1i_2}\dots \mathsf{r}_{i_ai}=1 && \textrm{\hspace{-6.5cm}for cycles with edges  $(ii_1),(i_1i_2)\dots,(i_ai)$}, 
\\[.6em]
\mathsf{r}_{ij}\mathsf{r}_{jk}\mathsf{t}_k=\mathsf{r}_{kj}\mathsf{r}_{ji}\mathsf{t}_i,\ \ \mathsf{t}_k\mathsf{r}_{ij}\mathsf{r}_{jk}=\mathsf{r}_{kj}\mathsf{r}_{ji}\mathsf{t}_i && \textrm{if $i<k$ and $m_{ik}=2$,}\\[.4em]
\mathsf{r}_{ij}\mathsf{r}_{jk}\mathsf{r}_{kl}\mathsf{t}_l=\mathsf{r}_{lk}\mathsf{r}_{kj}
\mathsf{r}_{ji}\mathsf{t}_i,\ \ \mathsf{t}_l\mathsf{r}_{ij}\mathsf{r}_{jk}\mathsf{r}_{kl}=\mathsf{r}_{lk}\mathsf{r}_{kj}
\mathsf{r}_{ji}\mathsf{t}_i && \textrm{if $i<l$ and $m_{il}=2$,}\\[.4em]
(\mathsf{r}_{ij}\mathsf{t}_j)^{\frac{m_{ij}}{2}}=(\mathsf{r}_{ji}\mathsf{t}_i)^{\frac{m_{ij}}{2}},\ \ (\mathsf{t}_j\mathsf{r}_{ij})^{\frac{m_{ij}}{2}}=(\mathsf{r}_{ji}\mathsf{t}_i)^{\frac{m_{ij}}{2}} && \text{if $m_{ij}>2$ and $m_{ij}$ is even,}\\[0.4em]
(\mathsf{r}_{ij}\mathsf{t}_j)^{\frac{m_{ij}-1}{2}}\mathsf{r}_{ij}=(\mathsf{r}_{ji}\mathsf{t}_i)^{\frac{m_{ij}-1}{2}},\ \ (\mathsf{r}_{ij}\mathsf{t}_j)^{\frac{m_{ij}+1}{2}}=\mathsf{t}_i(\mathsf{r}_{ji}\mathsf{t}_i)^{\frac{m_{ij}-1}{2}} && \text{if $m_{ij}>2$ and $m_{ij}$ is odd,}\\[0.6em]
\mathsf{r}_{ij}\mathsf{r}_{lm}=\mathsf{r}_{lm}\mathsf{r}_{ij} && \textrm{\hspace{-6.5cm}if $(ij)$ and $(lm)$ are not connected.}\end{array}\right.\end{equation}
\end{prop}
\emph{Proof.} The proof is similar to the proof of the Proposition \ref{pres-edge}. We skip the calculations and indicate below only the mutually inverse isomorphisms between the group 
$\beth$
and the group $\mathcal{B}^+(G)$ with the presentation of the Proposition \ref{pres-braid}:
\[\mathsf{r}_{ij}\mapsto \mathsf{R}'_i\mathsf{R}'^{-1}_j\ \ \ \text{and}\ \ \ \mathsf{t}_i\mapsto \mathsf{R}'_i\mathsf{R}_i,\ \, i=0,\dots,n-1,\]
and 
\[\mathsf{R}_0\mapsto \mathsf{t}_0,\ \,\mathsf{R}'_0\mapsto 1,\ \ \ \ \mathsf{R}_i\mapsto \mathsf{r}_{0i_1}\mathsf{r}_{i_1i_2}\dots \mathsf{r}_{i_ai}\mathsf{t}_i\ \,\text{and}\ \,\mathsf{R}'_i\mapsto \mathsf{r}_{ii_a}\dots \mathsf{r}_{i_2i_1}\mathsf{r}_{i_10},\ \,i=1,\dots,n-1,\]
where, for $i=1,\dots,n-1$, $(0,i_1,i_2,\dots,i_a,i)$ is a path in the graph $\mathcal{G}^c$ from the vertex $0$ to the vertex $i$. The second map is well-defined since the image of $\mathsf{R}_i$ (respectively, of $\mathsf{R}'_i$) does not depend on the chosen path, due to the first relation in (\ref{rel-braid2}).\hfill$\square$

\vskip .1cm
The isomorphism between the group generated by the elements $\mathsf{r}_{ij}$ and $\mathsf{t}_0,\dots,\mathsf{t}_{n-1}$ with the defining relations (\ref{rel-braid2}) and the subgroup $\mathcal{B}^+(G)$ of the braid group $\mathcal{B}(G)$ is given by:
\begin{equation}\label{iso-braid2}\mathsf{t}_i\mapsto g_i^2,\,\ i=0,\dots,n-1,\ \ \ \text{and}\ \ \ \mathsf{r}_{ij}\mapsto g_ig_j^{-1},\ \ \text{for all generators $\mathsf{r}_{ij}$}.\end{equation} 

\noindent
\textbf{Remarks.} \textbf{(i)} The action of the standard anti-automorphism $\tau$ (see remark \textbf{(i)} after the proof of the Proposition \ref{pres-braid}) on the generators of the presentation given by the Proposition \ref{prop-edge-br} is
\[\mathsf{t}_i\mapsto\mathsf{t}_i,\,\ i=0,\dots,n-1,\ \ \ \text{and}\ \ \ \mathsf{r}_{ij}\mapsto\mathsf{t}_j^{-1}\mathsf{r}_{ji}\mathsf{t}_i,\ \ \text{for all generators $\mathsf{r}_{ij}$}.\]

\textbf{(ii)} This remark is the analogue, for this presentation, of the remark \textbf{(ii)} after the proof of the Proposition \ref{pres-braid}. The alternating subgroup $G^+$ of the Coxeter group $G$ is the quotient of the group $\mathcal{B}^+(G)$, with the presentation (\ref{rel-braid2}), by the relations $\mathsf{t}_i=1$, $i=0,\dots,n-1$. Indeed, in this quotient, the relations (\ref{rel-braid2}) reduce immediately to the defining relations (\ref{code-edge}) of $G^+$.

\vskip .1cm
\textbf{(iii)} In the type A situation, with the same labeling of the Coxeter graph as in remark \textbf{(v)} after the proof of the Proposition \ref{pres-braid}, the presentation using edges of the Coxeter graph is the following. Set $\mathsf{r}_i:=\mathsf{r}_{i-1i}$, $i=1,\dots,n-1$. The group $\mathcal{B}^+(A_n)$ is isomorphic to the group generated by $\mathsf{r}_1,\dots,\mathsf{r}_{n-1}$ and $\mathsf{t}_0,\dots,\mathsf{t}_{n-1}$ with the defining relations:
\begin{equation}\label{braid2-A+}
\left\{\begin{array}{ll}
\mathsf{r}_i\mathsf{r}_{i+1}\mathsf{t}_{i+1}=\mathsf{r}_{i+1}^{-1}\mathsf{r}_i^{-1}\mathsf{t}_{i-1},\ \ 
\mathsf{t}_{i+1}\mathsf{r}_i\mathsf{r}_{i+1}=\mathsf{r}_{i+1}^{-1}\mathsf{r}_i^{-1}\mathsf{t}_{i-1} & \text{for $i=1,\dots,n-2$,}\\[0.3em]
\mathsf{r}_i\mathsf{r}_{i+1}\mathsf{r}_{i+2}\mathsf{t}_{i+2}=
\mathsf{r}_{i+2}^{-1}\mathsf{r}_{i+1}^{-1}\mathsf{r}_i^{-1}\mathsf{t}_{i-1},\ \ \mathsf{t}_{i+2}\mathsf{r}_i\mathsf{r}_{i+1}\mathsf{r}_{i+2}=
\mathsf{r}_{i+2}^{-1}\mathsf{r}_{i+1}^{-1}\mathsf{r}_i^{-1}\mathsf{t}_{i-1} & \text{for $i=1,\dots,n-3$,}\\[0.3em]
\mathsf{r}_i\mathsf{t}_i\mathsf{r}_i=\mathsf{r}_i^{-1}\mathsf{t}_i^{-1},\ \ (\mathsf{r}_i\mathsf{t}_{i})^2=\mathsf{t}_{i-1}\mathsf{r}_i^{-1}\mathsf{t}_{i-1} & \text{for $i=1,\dots,n-1$,}\\[0.4em]
\mathsf{r}_i\mathsf{r}_j=\mathsf{r}_j\mathsf{r}_i & \text{\hspace{-8cm}for $i,j=1,\dots,n-1$ such that $|i-j|>2$.}
\end{array}\right.
\end{equation} 
It is immediate to check (with the help of the isomorphism (\ref{iso-braid2})) that the following relations are satisfied
\[\mathsf{r}_i\mathsf{r}_j=\mathsf{r}_j\mathsf{r}_i,\ \ \mathsf{r}_i\mathsf{t}_j=\mathsf{t}_j\mathsf{r}_i,\ \ \mathsf{t}_i\mathsf{t}_j=\mathsf{t}_j\mathsf{t}_i\ \ \ \ \ \text{if $|i-j|>2$.}\]
So the presentation (\ref{braid2-A+}) of the chain of the groups $\mathcal{B}^+(A_n)$ is local and stationary,
in the sense of \cite{V}.

\section*{Appendix. \hspace{0.2cm} Coefficients in the defining relations of the alternating subalgebras}
\addcontentsline{toc}{section}{Appendix.$\ $ Coefficients in the defining relations of the alternating subalgebras\vspace{.1cm}}

The coefficients $a_k^{(m_{ij})}$
appearing in (\ref{mult-Hecke-f22}), (\ref{rel-bour}) and (\ref{rel-edge-a}) are easy to calculate for small $m_{ij}$. 
We define, in this Appendix, certain integers $\alpha_{k,l,l'}^{(m)}$ in terms of which the elements $a_k^{(m_{ij})}$, for any $m_{ij}$, can be expressed. We give the recurrent (in $m$) relations for
$\alpha_{k,l,l'}^{(m)}$ and find the generating function for $\alpha_{k,l,l'}^{(m)}$.
In the one-parameter situation, we recover a closed formula from \cite{Rat} for $a_k^{(m_{ij})}$.

\paragraph{Recursion.}
Let $f_1$ and $f_2$ be the generators of an algebra with the defining relations $f_1^2=f_2^2=1$. Define the elements $\alpha_{k,l,l'}^{(m)}\in\mathbb{Z}$, with $m,l,l'\in\mathbb{Z}_{\geq0}$ and $k\in\mathbb{Z}$, by
\begin{equation}\label{form-coeff}
\pr{(f_1+x)}{(f_2+y)}_m =\sum\limits_{l,l'\geq0}x^ly^{l'}\,\sum\limits_{k\in\mathbb{Z}}\alpha_{k,l,l'}^{(m)}\pr{f_1}{f_2}_k\ .\end{equation} 
By construction, for a given $m$, only a finite number of elements $\alpha_{k,l,l'}^{(m)}$ 
are non-zero; we have:
\[\alpha_{k,l,l'}^{(m)}\neq 0\ \ \Rightarrow\ \  
|k|\leq m-l-l',\ \ k\equiv m-l-l'\,({\rm mod}\, 2),\ \ l\leq\Res{(m+1)/2}\ \ \text{and}\ \ l'\leq\Res{m/2}\ .\]
Here $\Res{x}$ is the integer part of $x$.
The coefficients $a_k^{(m_{ij})}$, in terms of $\alpha_{k,l,l'}^{(m)}$,  read
\[a_k^{(m_{ij})}=\sum\limits_{l,l'\geq0}\beta_i^l\beta_j^{l'}\,(\alpha_{k,l,l'}^{(m_{ij})}-\alpha_{-k,l',l}^{(m_{ij})})\ .\]

\begin{lemm}\hspace{-.2cm}.\hspace{.2cm}
\label{lemm-coeff}
The elements $\alpha_{k,l,l'}^{(m)}$ satisfy the following initial condition and recursion:
\begin{eqnarray}
\label{rec-coeff1}
&&\alpha_{k,l,l'}^{(0)}=\delta^0_k\,\delta^0_l\,\delta^0_{l'}\ ,\\[0.7em]
\label{rec-coeff2}&&\alpha_{k,l,l'}^{(m+1)}=\alpha_{k-1,l',l}^{(m)}+\alpha_{-k,l',l-1}^{(m)}\ ,
\end{eqnarray}
where $\delta^i_j$ is the Kronecker delta.
\end{lemm}
\emph{Proof.} The initial condition 
(\ref{rec-coeff1}) is obviously verified. For the recurrence relation
(\ref{rec-coeff2}), one only has to notice that $\pr{(f_1+x)}{(f_2+y)}_{m+1} =(f_1+x)\pr{(f_2+y)}{(f_1+x)}_m$ and then use the induction hypothesis and $f_1\pr{f_2}{f_1}_k=\pr{f_1}{f_2}_{k+1}$ for any $k\in\mathbb{Z}$.\hfill$\square$

\paragraph{Generating function.}
Let
$C(t,u,v,s):=\sum\limits_{m,l,l'\geq0}\,\sum\limits_{k\in\mathbb{Z}}\alpha_{k,l,l'}^{(m)}\,t^mu^lv^{l'}s^k$. 
The formulas (\ref{rec-coeff1})--(\ref{rec-coeff2}) imply
\begin{equation}\label{rec-gen}
C(t,u,v,s)-1=tsC(t,v,u,s)+tuC(t,v,u,s^{-1})\ .
\end{equation}
Exchanging, in (\ref{rec-gen}),  $u$ and $v$, or replacing $s$ by $s^{-1}$, or doing both simultaneously, we obtain the following system of equations:
\[\left(\begin{array}{cccc}
1 & -ts & 0 & -tu\\
-ts & 1 & -tv & 0\\
0 & -tu & 1 & -ts^{-1}\\
-tv & 0 & -ts^{-1} & 1
        \end{array}\right) \left(\begin{array}{l}C(t,u,v,s)\\ C(t,v,u,s)\\ C(t,u,v,s^{-1})\\ C(t,v,u,s^{-1})\end{array}\right)=\left(\begin{array}{c}1\\ 1\\ 1\\ 1\end{array}\right)\ .\]
Inverting the matrix in the left hand side, we find the 
generating function
\begin{equation}\label{gen}
C(t,u,v,s)=\frac{1+t(u+s)+t^2(vs+us^{-1}-s^{-2}-uv)+t^3(1-u^2)(v-s^{-1})}{1-t^2(s^2+s^{-2}+2uv)+t^4(u^2-1)(v^2-1)}\ .
\end{equation}

\paragraph{One-parameter situation, $q_i=q_j$.}  
Now $\beta_i=\beta_j$ and the coefficients $a_k^{(m_{ij})}$ in terms of $\alpha_{k,L}^{(m)}$, read
\begin{equation}\label{eq-app1}a_k^{(m_{ij})}=\sum\limits_{L\geq0}\beta_i^L\,(\alpha_{k,L}^{(m_{ij})}-\alpha_{-k,L}^{(m_{ij})})\ \ ,\ \ \text{where}
\ \ \alpha_{k,L}^{(m)}:=\sum\limits_{l,l'\geq0\ :\ l+l'=L}\alpha_{k,l,l'}^{(m)}\ .
\end{equation}
Define $D(t,u,s):=\sum\limits_{m,L\geq0}\,\sum\limits_{k\in\mathbb{Z}}\alpha_{k,L}^{(m)}\,t^mu^Ls^k$.  
By definition, $D(t,u,s)=C(t,u,u,s)$. 
By (\ref{gen}),
\begin{equation}\label{gen-easy}
D(t,u,s)=\frac{1+t(u-s^{-1})}{1-t(s+s^{-1})+t^2(1-u^2)}=\frac{(1-ts^{-1})+tu}{(1-ts)(1-ts^{-1})}\ d(t,u,s)\ ,\end{equation}
where 
\begin{equation}\label{gen-easy0} d(t,u,s):= 
\left(1-{\displaystyle\frac{t^2u^2}{(1-ts)(1-ts^{-1})}}\right)^{-1}\ .
\end{equation}
The even in $u$ part of $D(t,u,s)$ is
\begin{equation}\label{gen-easy1} \frac{1}{1-ts}\ d(t,u,s)=
\sum\limits_{a,b,c\geq0}\frac{(a+b)!(a+c-1)!}{a!\,(a-1)!\,b!\,c!}\  t^{2a+b+c}\ u^{2a}\,s^{b-c}\ .
 \end{equation}
The odd in $u$ part of $D(t,u,s)$ is
\begin{equation}\label{gen-easy2}\frac{tu}{(1-ts)(1-ts^{-1})}
\ d(t,u,s)=
\sum\limits_{a,b,c\geq0}\frac{(a+b)!(a+c)!}{a!\,a!\,b!\,c!}\  t^{2a+1+b+c}\ u^{2a+1}\,s^{b-c}\ .
 \end{equation}
Thus $\alpha_{k,L}^{(m)}$ ($L,m\in\mathbb{Z}_{\geq0}$; $k\in\mathbb{Z}$) vanish unless $L\leq m$, $|k|\leq m-L$ and $m+k\equiv L\,({\rm mod}\, 2)$, and then
\begin{equation}\label{11}
\left\{\begin{array}{lcl}
\alpha_{k,\,L}^{(m)}=\left(\begin{array}{c}(m+k)/2\\
(m+k-L)/2\end{array}\right) \left(\begin{array}{c}(m-k-2)/2\\
(m-k-L)/2\end{array}\right)&,&  L\equiv 0\,({\rm mod}\, 2)\ ,
\\[1.5em]
\alpha_{k,\,L}^{(m)}=\left(\begin{array}{c}(m+k-1)/2\\
(m+k-L)/2\end{array}\right) \left(\begin{array}{c}(m-k-1)/2\\
(m-k-L)/2\end{array}\right) \hspace{0.2cm}&,&  L\equiv 1\,({\rm mod}\, 2)\ .
\end{array}\right.\end{equation}
{}For $k\geq 0$, $\alpha_{-k,\,L}^{(m)}=\alpha_{k,\,L}^{(m)}$ if $L\equiv 1\,({\rm mod}\, 2)$, and $\alpha_{-k,\,L}^{(m)}=\displaystyle\frac{m-k}{m+k}\alpha_{k,\,L}^{(m)}$ if $L\equiv 0\,({\rm mod}\, 2)$.
Substituting into (\ref{eq-app1}), one finds
\vskip -.5cm
\begin{equation}\label{for-a}
a_k^{(m_{ij})}=\sum\limits_{p=0}^{\Res{(m-1)/2}}\beta_i^{2p}\ \frac{2k}{m_{ij}+k}\ \alpha_{k,\,2p}^{(m_{ij})}\ \ \ \text{for $k=1,\dots,m_{ij}$.}
\end{equation}
The formulas (\ref{11}) and (\ref{for-a}) appear in \cite{Rat}.

\paragraph{Examples.}
We write explicitly the relation (\ref{rel-edge-a}) for $m_{ij}\leq6$ (this is all what is needed for the finite Coxeter groups other than the dihedral groups ${\sf D}_n$ with $n>6$) in the multiparameter setting:

\vskip .05cm
\noindent $\ \ \bullet$ for $m_{ij}=2$, $y_{ij}^2=1$ ;

\vskip .05cm
\noindent $\ \ \bullet$ for $m_{ij}=3$, $y_{ij}^3=\beta_i^2(y_{ij}-y_{ij}^2)+1$ ;

\vskip .05cm
\noindent $\ \ \bullet$ for $m_{ij}=4$, $y_{ij}^4=2\beta_i\beta_j(y_{ij}-y_{ij}^3)+1$ ;

\vskip .05cm
\noindent $\ \ \bullet$ for $m_{ij}=5$, $y_{ij}^5=3\beta_i^2(y_{ij}-y_{ij}^4) +(\beta_i^4+\beta_i^2)(y_{ij}^2-y_{ij}^3)+1$ ;

\vskip .05cm
\noindent $\ \ \bullet$ for $m_{ij}=6$, $y_{ij}^6=4\beta_i\beta_j(y_{ij}-y_{ij}^5) +(3\beta_i^2\beta_j^2+\beta_i^2+\beta_j^2)(y_{ij}^2-y_{ij}^4)+1$ .

\vskip .1cm
\noindent The one-parameter situation is obtained when one sets $\beta_i=\beta_j$.

\paragraph{Acknowledgment.} We thank I. Marin for useful discussions.

\end{document}